\documentclass[a4paper,12pt]{article}
\usepackage[english]{babel}   
\usepackage{indentfirst}
\usepackage[T2A]{fontenc}
\usepackage[utf8]{inputenc}

\renewcommand{\epsilon}{\ensuremath{\varepsilon}}
\renewcommand{\phi}{\ensuremath{\varphi}}
\renewcommand{\kappa}{\ensuremath{\varkappa}}
\renewcommand{\le}{\ensuremath{\leqslant}}

\renewcommand{\ge}{\ensuremath{\geqslant}}

\usepackage{amsmath,amsfonts,amssymb,amsthm,mathtools} 
\usepackage{icomma} 

\usepackage{graphicx}  
\graphicspath{{images/}{images2/}}  
\usepackage{geometry} 
\usepackage{setspace} 
\usepackage{lastpage} 
\usepackage{soul} 
  \usepackage{hyperref}
\usepackage[usenames,dvipsnames,svgnames,table,rgb]{xcolor}
\begin{document}

\begin{center} {\huge Legendre's conjecture. Theorem on existence
\begin{spacing}{1.3}\end{spacing}
of a prime number between} \LARGE{${\rm \;}m^{2}$} \huge{and} \LARGE{$(m+1)^{2}$}
\end{center}

\textbf{}

\begin{center} {\large Garipov Ilshat Ilsurovich }\end{center}
\begin{spacing}{-1.5}\end{spacing}
\begin{center} \textit {\normalsize Russia, Republic of Tatarstan, Naberezhnye Chelny}\end{center}
\begin{spacing}{-1.5}\end{spacing}
\begin{center} \textit {\normalsize e-mail: mathsciencegaripovii@gmail.com.}\end{center}
\renewcommand{\abstractname}{Abstract}
\begin{spacing}{0.5}\end{spacing}
\begin{abstract}
{\small In scientific paper, we will show a proof of Legendre's conjecture based on a scheme for finding elements of «active» set ${\rm H} _{m^{4} }$ and «critical» element $C_{m^{4}}$ for number $m^{4}$ at each $m \ge 3$.}

\end{abstract}

\begin{center}
\textbf{\large List of symbols}
\end{center}
\begin{spacing}{0.8}\end{spacing}

${\rm \mathbb{N}} $ -- set of all natural numbers.

${\rm \mathbb{P}} $ -- set of all prime numbers.

${\rm \mathbb{R}}$ -- set of all real numbers.

$T$ -- matrix comprising all defining and all not defining elements.

$D(b)$ -- $T$-matrix upper defining element of number $b$.

$D_{T} $ -- set of all defining elements of $T$-matrix.

$nD_{T} $ -- set of all not defining elements of $T$-matrix.

${\rm M} _{T} $ -- set of all leading elements of $T$-matrix.

$D_{T_{k} } $ -- set of all defining elements in $k$-row ($k>1$) of $T$-matrix.

$\pi (x)$ -- function counting the number of prime numbers less than or equal to $x\in {\rm \mathbb{R}}$.

$\pi _{{\rm M} _{T} } (x)$ -- function counting the number of $T$-matrix leading elements less than or equal to $x\in {\rm \mathbb{R}}$.

$\# _{k} (a)$ -- number of element $a$ in $k$-row of $T$-matrix.

${\rm H} _{m^{4} } $ -- «active» set for numbers $m^{4}$.

$C_{m^{4} } $ -- «critical» element for numbers $m^{4}$.

$q_{m}$ -- number of prime numbers between $m^{2}$ and $(m+1)^{2} $.

\newpage

\begin{center}
\textbf{\large Introduction}
\end{center}
\begin{spacing}{0.8}\end{spacing}

In paper [1], we explored a $T$-matrix: examined the simplest properties of $T$-matrix, proved the basic theorems about elements of $T$-matrix, found her applications in number theory. In paper [2], we established the connection between $T$-matrix and Legendre's conjecture, got «Weak» and «Strong» conjectures.

In this paper, we will present one proof of Legendre's conjecture using a scheme for finding ${\rm H} _{m^{4} }, C_{m^{4}}$ for number $m^{4}$ at each $m \ge 3$.

\begin{center}
\textbf{\large 1. About elements of $T$-matrix}
\end{center}
\begin{spacing}{0.7}\end{spacing}

We will first introduce the support constructions, definitions, propositions from [1], [2].

Construct a matrix $T\equiv \left(a(k;n)\right)_{\infty \times \infty } $, where $a(k;n)$ is a $T$-matrix element located in $k$-th row, $n$-th column and defined as follows:
\begin{spacing}{0.8}\end{spacing}
\[a(k;n)\equiv p(k)\cdot \left(5+2\cdot \left\lfloor \frac{n}{2} \right\rfloor +4\cdot \left\lfloor \frac{n-1}{2} \right\rfloor \right),\] 
\noindent where $p(k)$ is the $k$-th element of sequence $(p(k))_{k=1}^{\infty }$ of prime numbers: 
\begin{spacing}{0.5}\end{spacing}
\begin{equation} \label{(1)} 
p(k)\equiv p_{k+2} ,   
\end{equation}\begin{spacing}{1.2}\end{spacing}
\noindent where $p_{i}$ is the $i$-th prime number in sequence of all prime numbers.

Let $(f(n))_{n=1}^{\infty } $ is a numerical sequence, where a common member $f(n)$ is defined as follows:
\begin{spacing}{0.8}\end{spacing}\[f(n)\equiv 3n+\frac{3-(-1)^{n} }{2} .\]\begin{spacing}{1.6}\end{spacing}
THEOREM 1.1. 
\[\left(\forall k,n\in {\rm \mathbb{N}} \right)\left(a(k;n)=p(k)\cdot f(n)\right).\]

DEFINITION 1.1. An element $a(k;n)$ of matrix $T$ is called defining $(a(k;n) \in D_{T})$ if

1)  $a(k;n)$ is not divisible by 5;

2) $a(k;n)$ can be expressed as a product of some two prime numbers, that is
\begin{spacing}{0.4}\end{spacing}\begin{equation} \label{(2)} 
5\not |a(k;n){\rm \; \; }\wedge {\rm \; \; }(\exists u ,v \in {\rm \mathbb{P}} )(a(k;n)=u \cdot v ).  
\end{equation}\begin{spacing}{1.1}\end{spacing}

DEFINITION 1.2. An element $a(k;n)$ of matrix $T$ is called not defining $(a(k;n) \in nD_{T})$ if he does not satisfy condition \eqref{(2)}.

DEFINITION 1.3. An element $a(k;n)$ of matrix $T$ is called leading $(a(k;n) \in {\rm M} _{T})$ if
\begin{spacing}{0.5}\end{spacing}
\[a(k;n)=p^{2} (k).\]

DEFINITION 1.4. A $T$-matrix is called matrix comprising all defining and not defining elements.

LEMMA 1.2. $(f(n))_{n=1}^{\infty } $ is a sequence of all numbers of the form $6h\pm 1$:
\begin{spacing}{0.5}\end{spacing}
\[5;{\rm \; }7;{\rm \; }11;{\rm \; }13;{\rm \; }17;{\rm \; }19;{\rm \; }23;{\rm \; }25;{\rm \; }...{\rm \; };{\rm \; }6h-1;{\rm \; 6}h+{\rm 1;...}.\]
\begin{spacing}{1.1}\end{spacing}

DEFINITION 1.5. A $T$-matrix defining element $D(b)$ is called an upper defining element of number $b\in {\rm \mathbb{R}}: b \ge 49$, if 
\begin{spacing}{0.5}\end{spacing}
\[D(b)=\min_{\scriptstyle \substack{a(k_{1}; n)\in D_{T} \\  a(k_{1}; n)>b \\  n\in {\rm \mathbb{N}}}}a(k_{1}; n), \]      
 \begin{spacing}{1.3}\end{spacing}

\noindent where $k_{1} $  is defined by condition
\begin{spacing}{0.3}\end{spacing}
\[ p^{2} (k_{1})=\max_{\scriptstyle \substack{p^{2} (k)\le b \\  k>1 }}{\rm \;}p^{2}(k).\]

THEOREM 1.3 (about the «transition down» of $T$-matrix defining element). 
\begin{spacing}{0.5}\end{spacing}\[\left(\forall k;n\in {\rm \mathbb{N}} \right)\left(p^{2} (k)<a(k;n){\rm \; \; }\wedge {\rm \; \; }a(k;n)\in D_{T} {\rm \;}\Rightarrow \right. \] 
\[\left. \Rightarrow \left(\exists ! {\rm \; }j\in {\rm \mathbb{N}}\right)\left(k<j{\rm \;}\wedge {\rm \;}a(k;n)<p^{2} (j){\rm \;}\wedge {\rm \;}a(j;\# _{k} (p^{2} (k)))=a(k;n){\rm \;}\wedge {\rm \;}a(j;n)=p^{2} (j)\right)\right).\]

\begin{center}
\textbf{\large 2. About scheme for getting ${\rm H} _{m^{4}}, C_{m^{4}}$ for number} $m^{4}{\rm \; } (m \ge 3)$ 
\end{center}
\begin{spacing}{0.7}\end{spacing}

Let any natural number $m\ge 3$ is chosen, and the leading element $p^{2} (k_{1})$ of $T$-matrix is found by condition:
\begin{spacing}{0.5}\end{spacing}
\begin{equation} \label{(3)} 
p^{2}(k_{1})=\max_{\scriptstyle \substack {p^{2}(k)<m^{4} \\ k>1}}{\rm \; } p^{2}(k), \text{where }k_{1} \equiv k_{1}(m).
\end{equation} 

THEOREM 2.1.
\begin{spacing}{0.1}\end{spacing} \begin{equation} \label{(4)}
\left(\forall x\in {\rm \mathbb{R}}\right)\left(x\ge 3{\rm \;}\Rightarrow {\rm \;}\pi (x)=\pi _{{\rm M} _{T} } \left(x^{2} \right)+2\right) (\text{see }[1]).         
\end{equation}\begin{spacing}{1}\end{spacing}

Now find a $k_{1}$-row of $T$-matrix.
\[\pi _{{\rm M} _{T} } (p^{2}(k_{1})){\rm \; \; }\mathop{=}\limits^{\eqref{(3)}}{\rm \; \; }\pi _{{\rm M} _{T} } \left (\max_{\scriptstyle \substack {p^{2}(k)<m^{4} \\ k>1}}{\rm \; } p^{2}(k)\right)=\max_{\scriptstyle \substack {p^{2}(k)<m^{4} \\ k>1}}{\rm \; } \pi _{{\rm M} _{T} } (p^{2}(k)){\rm \; \; }\mathop{=}\limits^{\eqref{(4)}}\]
\[=\max_{\scriptstyle \substack {p^{2}(k)<m^{4} \\ k>1}}{\rm \; } (\pi (p(k))-2){\rm \; \; }\mathop{=}\limits^{\eqref{(1)}}{\rm \;}\max_{\scriptstyle \substack {p^{2}(k)<m^{4} \\ k>1}}{\rm \; } (\pi (p_{k+2})-2)=\]

\[=\max_{\scriptstyle \substack {p^{2}(k)<m^{4} \\ k>1}}{\rm \; } ((k+2)-2)=\max_{\scriptstyle \substack {p^{2}(k)<m^{4} \\ k>1}}{\rm \; } k.\]
\[\pi _{{\rm M} _{T} } (p^{2}(k_{1})) {\rm \; \; }\mathop{=} \limits^{\eqref{(4)}}{\rm \; \; }\pi (p(k_{1}))-2 {\rm \; \; }\mathop{=} \limits^{\eqref{(1)}}{\rm \; \; }\pi (p_{k_{1}+2})-2=(k_{1}+2)-2=k_{1}.\]

It follows that
\begin{equation} \label{(5)}
k_{1}(m)=\max_{\scriptstyle \substack {p^{2}(k)<m^{4} \\ k>1}}{\rm \; } k.
\end{equation}

In paper [2], the definition of «active» set ${\rm H} _{(m-1)^{4} ,{\rm \; }m^{4}}$ and definition of «critical» element $C_{(m-1)^{4} ,{\rm \;}m^{4}}$ for numbers $(m-1)^{4} ,{\rm \;}m^{4}$ $(m \ge 3)$ were introduced under the assumption that Legendre's conjecture is true for number $m-1$. By removing this assumption, we will formulate 3 definitions.

DEFINITION 2.1. A set  ${\rm H} _{m^{4}}\subset$  $D_{T_{k_{1}(m) } } $ of all defining elements  $a(k_{1};n_{i})>p^{2} (k_{1})$ (in $k_{1}$-row of $T$-matrix), $i=\overline{1;q_{m}}$ , in which 

\[\left[\begin{array}{l} {a(k_{1};n_{i})<m^{4} <p^{2} (j_{i})<(m+1)^{4}} \\ {m^{4} <a(k_{1} ;n_{i})<p^{2} (j_{i})<(m+1)^{4}} \end{array}\right. \text{,}\]

\noindent where $n_{1}<n_{2}<...<n_{q_{m}}, j_{i}=k_{1}+i$; is called an «active» set for number $m^{4}$ ($m\ge 3$).

DEFINITION 2.2. A defining element  $C_{m^{4}} \equiv a(k_{1} ;n_{q_{m}+1})\notin {\rm H}_{m^{4} } $, next to a defining element $a(k_{1} ;n_{q_{m}} )\in {\rm H} _{m^{4} } {\rm \;}(q_{m}>0)$, is called a «critical» element for number $m^{4}$ ($m\ge 3$).

Take $C_{m^{4}} \equiv D(p^{2}(k_{1}))$ if $q_{m}=0$ for some $m \ge 3$. In this case, ${\rm H} _{m^{4} }=\varnothing$.
\begin{spacing}{1.2}\end{spacing}

DEFINITION 2.3. A «transition» of the defining element $a(k_{1} ;n_{i})>p^{2}(k_{1})$ from $k_{1}$-row to $j_{i}$-row ($j_{i}>k_{1}$) of $T$-matrix with some $i\in {\rm \mathbb{N}}$ is called successful if $a(k_{1};n_{i})\in {\rm H} _{(m-1)^{4} ,{\rm \; }m^{4} } $. Otherwise, that is if $a(k_{1} ;n_{i})\notin {\rm H}_{(m-1)^{4} ,{\rm \; }m^{4}}$, this «transition» is called unsuccessful.

\textbf{Scheme №1 for getting ${\rm H} _{m^{4}}, C_{m^{4}}$ for number $m^{4}$ $(m \ge 3, k_{1} \equiv k_{1}(m))$}.
\begin{spacing}{-0.3}\end{spacing}
\[p^{2}(k_{1}(m)){\rm \; } \mathop{ \rightarrow}\limits^{k_{1} \rightarrow k_{1}}{\rm \; } a(k_{1} ;n_{1}){\rm \; } \mathop{ \rightarrow}\limits^{k_{1} \rightarrow  k_{1}+1} {\rm \; }a(k_{1}+1;\# _{k_{1}} (p^{2} (k_{1}))){\rm \; }\mathop{ \rightarrow}\limits^{k_{1}+1 \rightarrow k_{1}}{\rm \; } a(k_{1} ;n_{2}){\rm \; }\mathop{ \rightarrow}\limits^{k_{1} \rightarrow k_{1}+2}\]
\begin{spacing}{-0.3}\end{spacing}
\[\mathop{ \rightarrow}\limits^{k_{1} \rightarrow k_{1}+2}{\rm \; } a(k_{1}+2;\# _{k_{1}} (p^{2} (k_{1}))){\rm \; }\mathop{ \rightarrow}\limits^{k_{1}+2 \rightarrow k_{1}}{\rm \; }...{\rm \; }\mathop{ \rightarrow}\limits^{k_{1}+q_{m}-1 \rightarrow k_{1}}{\rm \; } a(k_{1} ;n_{q_{m}}){\rm \; }\mathop{ \rightarrow}\limits^{ k_{1} \rightarrow k_{1}+q_{m}}\]
\[\mathop{ \rightarrow}\limits^{ k_{1} \rightarrow k_{1}+q_{m}} {\rm \; }a(k_{1}+q_{m};\# _{k_{1}} (p^{2} (k_{1}))){\rm \; }\mathop{ \rightarrow}\limits^{ k_{1}+q_{m} \rightarrow k_{1} }{\rm \; }a(k_{1} ;n_{q_{m}+1}){\rm \; }\mathop{ \rightarrow}\limits^{ k_{1}  \rightarrow k_{1}+q_{m}+1}\]
\[\mathop{ \rightarrow}\limits^{ k_{1} \rightarrow k_{1}+q_{m}+1}{\rm \; }a(k_{1}+q_{m}+1;\# _{k_{1}} (p^{2} (k_{1}))){\rm \; }\mathop{ \rightarrow}\limits^{ k_{1}+{q_{m}+1} \rightarrow k_{1}+q_{m}}{\rm \; }p^{2}(k_{1}(m+1)).\]

Here $a \mathop{ \rightarrow}\limits^{ i \rightarrow j} b$ indicates the transition from element $a \in D_{T_{i}}$ in $i$-row to element $b \in D_{T_{j}}$ in $j$-row of $T$-matrix.

$a(k_{1};n_{i+1})$ is a defining element, next to a defining element $a(k_{1};n_{i})$, $i \in \mathbb{N}$.

From the above Scheme №1, we get:

1) ${\rm H} _{m^{4}}={\rm \; \{ }a(k_{1} ;n_{1} ),{\rm .}..,{\rm \; }a(k_{1} ;n_{q_{m}} )\}$ , where $p^{2}(k_{1})<a(k_{1} ;n_{1} )<...<a(k_{1} ;n_{q_{m}} )$.

2) $C_{m^{4}}=a(k_{1} ;n_{q_{m}+1})$, where $C_{m^{4}}>a(k_{1} ;n_{q_{m}} )$.

3) $p^{2}(k_{1}(m+1)) = a(k_{1}+{q_m};n_{q_m})$, where 

\begin{spacing}{0.3}\end{spacing}
\begin{equation} \label{(6)} 
k_{1}(m+1)=k_{1}(m)+q_{m}.
\end{equation}
\begin{spacing}{1}\end{spacing}

\textbf{Description of Scheme №1. }Find a leading element $p^{2}(k_{1}(m))$ in $k_{1}$-row of $T$-matrix for a given number $m \ge 3$. It should be reminded that $p^{2}(k_{1}(m))$ is maximal leading element that is less than $m^{4}$. In $k_{1}$-row of $T$-matrix, move from leading element $p^{2}(k_{1}(m))$ to defining element $a(k_{1} ;n_{1})$.

Move from $a(k_{1} ;n_{1})$ in $k_{1}$-row of $T$-matrix to defining element $a(k_{1}+1;\# _{k_{1}} (p^{2} (k_{1})))$ in $k_{1}+1$-row of $T$-matrix. By Theorem 1.3,
\begin{spacing}{0.5}\end{spacing}
\[a(k_{1}+1;\# _{k_{1}} (p^{2} (k_{1})))=a(k_{1} ;n_{1}).\]

The transition of defining element $a(k_{1} ;n_{1})>p^{2}(k_{1})$ from $k_{1}$-row to $k_{1}+1$-row of $T$-matrix is successful. Therefore we will move from defining element $a(k_{1}+1;\# _{k_{1}} (p^{2} (k_{1})))$ in $k_{1}+1$-row of $T$-matrix to defining element $a(k_{1};n_{2})$ in $k_{1}$-row of $T$-matrix.

Move from $a(k_{1} ;n_{2})$ in $k_{1}$-row of $T$-matrix to defining element $a(k_{1}+2;\# _{k_{1}} (p^{2} (k_{1})))$ in $k_{1}+2$-row of $T$-matrix. By Theorem 1.3,
\begin{spacing}{0.5}\end{spacing}
\[a(k_{1}+2;\# _{k_{1}} (p^{2} (k_{1})))=a(k_{1} ;n_{2}).\]

The transition of defining element $a(k_{1} ;n_{2})>p^{2}(k_{1})$ from $k_{1}$-row to $k_{1}+2$-row of $T$-matrix is successful. Therefore we will move from defining element $a(k_{1}+2;\# _{k_{1}} (p^{2} (k_{1})))$ in $k_{1}+2$-row of $T$-matrix to defining element $a(k_{1};n_{3})$ in $k_{1}$-row of $T$-matrix. And so on to defining element $a(k_{1} ;n_{q_{m}})$ in $k_{1}$-row of $T$-matrix.

Move from $a(k_{1} ;n_{q_{m}})$ in $k_{1}$-row of $T$-matrix to defining element $a(k_{1}+q_{m};\# _{k_{1}} (p^{2} (k_{1})))$ in $k_{1}+q_{m}$-row of $T$-matrix. By Theorem 1.3,
\begin{spacing}{0.5}\end{spacing}
\[a(k_{1}+q_{m};\# _{k_{1}} (p^{2} (k_{1})))=a(k_{1} ;n_{q_{m}}).\]

The transition of defining element $a(k_{1} ;n_{q_{m}})>p^{2}(k_{1})$ from $k_{1}$-row to $k_{1}+q_{m}$-row of $T$-matrix is successful. Therefore we will move from defining element  $a(k_{1}+q_{m};\# _{k_{1}} (p^{2} (k_{1})))$ in  $k_{1}+q_{m}$-row of $T$-matrix to defining element $a(k_{1} ;n_{q_{m}+1})$ in $k_{1}$-row of $T$-matrix.

Move from $a(k_{1} ;n_{q_{m}+1})$ in $k_{1}$-row of $T$-matrix to defining element $a(k_{1}+q_{m}+1;\# _{k_{1}} (p^{2} (k_{1})))$ in $k_{1}+q_{m}+1$-row of $T$-matrix. By Theorem 1.3,
\begin{spacing}{0.5}\end{spacing}
\[a(k_{1}+q_{m}+1;\# _{k_{1}} (p^{2} (k_{1})))=a(k_{1} ;n_{q_{m}+1}).\]

The transition of defining element $a(k_{1} ;n_{q_{m}+1})>p^{2}(k_{1})$ from $k_{1}$-row to $k_{1}+q_{m}+1$-row of $T$-matrix is unsuccessful. In this case, we will move from $a(k_{1}+q_{m}+1;\# _{k_{1}} (p^{2} (k_{1})))$ in $k_{1}+q_{m}+1$-row of $T$-matrix to leading element  $p^{2}(k_{1}(m+1))$ in $k_{1}+q_{m}$-row of $T$-matrix. In turn, $p^{2}(k_{1}(m+1))$ is maximal leading element that is less than $(m+1)^{4}$.

\begin{center}
\textbf{\large 3. Proof of Legendre's conjecture}
\end{center}

\begin{spacing}{0.7}\end{spacing}
THEOREM 3.1. Set ${\rm M}_{T}$ is infinite (see [1]).

Note that properties and propositions about «active» set ${\rm H} _{(m-1)^{4} ,{\rm \; }m^{4}}$ for numbers $(m-1)^{4}$, $m^{4}$ (see [2]) will also be true for «active» set ${\rm H} _{m^{4}}$ for number $m^{4}$ $(m\ge 3)$.

Definition 2.1 makes it clear that $\left(\forall m\ge 3\right)\left(\left|{\rm H} _{m^{4}} \right|=q_{m}\right)$. Therefore, set ${\rm H} _{m^{4}}$ is finite.

Let's define the set $D_{T_{k_{1}(m)}}^{*}$ as follows:
\[D_{T_{k_{1}(m)}}^{*}  \equiv \left\{ a \in D_{T_{k_{1}(m)}} {\rm \;} | {\rm \; \;} p^{2}(k_{1}(m))<a \le C_{m^{4}}\right\}.\]

It follows from Definition 2.1 that 
\[{\rm H} _{m^{4}}=\left\{ a \in D_{T_{k_{1}(m)}} {\rm \;} | {\rm \; \;} p^{2}(k_{1}(m))<a< C_{m^{4}}\right\}.\]

In turn,
\[\overline{{\rm H}} _{m^{4}} \equiv D_{T_{k_{1}(m)}}^{*} \backslash {\rm \;}{\rm H} _{m^{4}}=\left\{ C_{m^{4}}\right\}.\]

THEOREM 3.2 (theorem on existence of a prime number between $m^{2}$ and $(m+1)^{2}$). 
\[\left(\forall m\in {\rm \mathbb{N}} \right)\left(\exists p\in {\rm \mathbb{P}} \right)\left(p \in \left(m^{2}; (m+1)^{2}\right) \right).\]

PROOF.  The Legendre's conjecture is true for $m \in \{1;2\}$. For example, the prime number $p=2$ lies between $1^{2}$ and $2^{2}$, the prime number $p=5$ lies between $2^{2}$ and $3^{2}$. 

Further, we consider the Legendre's conjecture for $m \ge 3$.

Suppose otherwise: $\left(\exists m \ge 3\right)\left(\forall p\in {\rm \mathbb{P}} \right)\left(p \not \in \left(m^{2}; (m+1)^{2}\right) \right).$ This is equivalent to the next:
\begin{spacing}{0.5}\end{spacing}
\[\left(\exists m \ge 3\right)\left(\forall q\in {\rm M}_{T} \right)\left(q \not \in \left(m^{4}; (m+1)^{4}\right) \right).\]

Let's assume we found a number $m_{0} \ge 3$ such that
\begin{spacing}{0.5}\end{spacing}
\[\left(\forall q\in {\rm M}_{T} \right)\left(q \not \in \left(m_{0}^{4}; (m_{0}+1)^{4}\right) \right).\]

Using \eqref{(5)}, we find  $k_{1}(m_{0})$-row of $T$-matrix:
\begin{spacing}{0.5}\end{spacing}
\[k_{1}(m_{0})=\max_{\scriptstyle \substack {p^{2}(k)<m_{0}^{4} \\ k>1}}{\rm \; } k.\]

Then,
\begin{spacing}{0.5}\end{spacing}
\[p^{2}(k_{1}(m_{0}))=\max_{\scriptstyle \substack {p^{2}(k)<m_{0}^{4} \\ k>1}}{\rm \; } p^{2}(k).\]

Given Theorem 3.1,
\begin{spacing}{0.5}\end{spacing}
\[\left(\exists s \in \mathbb{N} \right)\left(p^{2}(k_{1}(m_{0})+1) \not \in \left(m_{0}^{4}; (m_{0}+s)^{4} \right) {\rm \; }\wedge \right.\]
\[ \left. \wedge {\rm \; \; } p^{2}(k_{1}(m_{0})+1) \in \left((m_{0}+s)^{4}; (m_{0}+s+1)^{4} \right) \right).\]

Fix such found number $s$. Let $k_{1} \equiv k_{1}(m_{0})$. 

Using Scheme №1, construct scheme for $m_{0}^{4}$, starting with leading element $p^{2}(k_{1}(m_{0}))$:
\[p^{2}(k_{1}(m_{0})){\rm \; \; } \mathop{ \rightarrow}\limits^{k_{1} \rightarrow k_{1}}{\rm \; \; }a(k_{1};n_{1}){\rm \; \; } \mathop{ \rightarrow}\limits^{k_{1} \rightarrow k_{1}+1}{\rm \; \; }a(k_{1}+1;\# _{k_{1}} (p^{2} (k_{1}))){\rm \; \; }\mathop{ \rightarrow}\limits^{ k_{1}+1 \rightarrow k_{1}}{\rm \; \; }p^{2}(k_{1}(m_{0}+1)).\]
\[k_{1}(m_{0}+1){\rm \;}\mathop{=}\limits^{\eqref{(6)}}{\rm \;}k_{1}(m_{0})+q_{m_{0}}=k_{1}(m_{0})+0=k_{1}(m_{0}){\rm \; \; }\wedge\]
\[\wedge {\rm \; \;} a(k_{1}+1;\# _{k_{1}} (p^{2} (k_{1})))=a(k_{1};n_{1})=D(p^{2}(k_{1})){\rm \; \;} \Rightarrow\]
\begin{equation} \label{(7)} 
\Rightarrow{\rm \; \; }p^{2}(k_{1}(m_{0})){\rm \; \; } \mathop{ \rightarrow}\limits^{k_{1} \rightarrow k_{1}}{\rm \; \; }D(p^{2}(k_{1})){\rm \; \; }\mathop{ \rightarrow}\limits^{ k_{1} \rightarrow k_{1}+1}{\rm \; \; }D(p^{2}(k_{1})){\rm \; \; }\mathop{ \rightarrow}\limits^{ k_{1}+1 \rightarrow k_{1}}{\rm \; \; }p^{2}(k_{1}(m_{0})).
\end{equation}

Scheme \eqref{(7)} will be the same for numbers $m_{0}^{4},...,(m_{0}+s-1)^{4}$ because 
\begin{spacing}{0.5}\end{spacing}
\begin{equation} \label{(8)} 
k_{1}(m_{0})=k_{1}(m_{0}+1){\rm \;}\mathop{=}\limits^{\eqref{(6)}}{\rm \;}...{\rm \;}\mathop{=}\limits^{\eqref{(6)}}{\rm \;}k_{1}(m_{0}+s).
\end{equation}

With such found number $s$ we get from scheme \eqref{(7)} and from equalities \eqref{(8)}:

1) ${\rm H} _{(m_{0}+i)^{4}}=\varnothing, i=\overline{0;s-1}.$

2) $C_{(m_{0}+i)^{4}}=D(p^{2}(k_1(m_0)))$,where $C_{(m_{0}+i)^{4}}>p^{2}(k_1(m_0)), i=\overline{0;s-1}.$

3) $p^{2}(k_{1}(m_{0}))=p^{2}(k_{1}(m_{0}+1))=...=p^{2}(k_{1}(m_{0}+s)).$

Introduce, $m \equiv m_{0}+s-1$. It can be considered that $k_{1} \equiv k_{1}(m)$. In this case, there is a following scheme for number $m^4$:
\begin{spacing}{0.5}\end{spacing}
\[p^{2}(k_{1}(m)){\rm \; \; } \mathop{ \rightarrow}\limits^{k_{1} \rightarrow k_{1}}{\rm \; \; }D(p^{2}(k_{1})){\rm \; \; }\mathop{ \rightarrow}\limits^{ k_{1} \rightarrow k_{1}+1}{\rm \; \; }D(p^{2}(k_{1})){\rm \; \; }\mathop{ \rightarrow}\limits^{ k_{1}+1 \rightarrow k_{1}}{\rm \; \; }p^{2}(k_{1}(m+1)). \]
\begin{equation} \label{(9)} 
k_{1}(m+1){\rm \;}\mathop{=}\limits^{\eqref{(6)}}{\rm \;}k_{1}(m)+q_{m}=k_{1}(m)+0=k_{1}(m) {\rm \;}\Rightarrow{\rm \;} k_{1}(m+1)=k_{1}(m).
\end{equation}
\[D_{T_{k_{1}(m)}}^{*}=\left\{ a \in D_{T_{k_{1}(m)}} {\rm \;} | {\rm \; \;} p^{2}(k_{1}(m))<a \le C_{m^{4}}\right\},\]
\[C_{m^{4}}=D(p^{2}(k_1(m)))=p(k_1(m)) \cdot  p(k_1(m)+1).\]
\begin{equation} \label{(10)} 
{\rm H} _{m^{4}}=\varnothing, \overline{{\rm H}} _{m^{4}}=\left\{ C_{m^{4}}\right\} {\rm \; \;}\Rightarrow {\rm \; \;}D_{T_{k_{1}(m)}}^{*}={\rm H} _{m^{4}} \cup \overline{{\rm H}} _{m^{4}}=\overline{{\rm H}} _{m^{4}}.
\end{equation}

Using Scheme №1 again, construct scheme for $(m+1)^{4}$, starting with leading element 
\begin{spacing}{1.5}\end{spacing}
\noindent $p^{2}(k_{1}(m+1))\mathop{=}\limits^{\eqref{(9)}}p^{2}(k_{1}(m))$.
\[p^{2}(k_{1}(m+1)){\rm \; } \mathop{ \rightarrow}\limits^{k_{1} \rightarrow k_{1}}{\rm \; } a(k_{1} ;n_{1}){\rm \; } \mathop{ \rightarrow}\limits^{k_{1} \rightarrow  k_{1}+1} {\rm \; }a(k_{1}+1;\# _{k_{1}} (p^{2} (k_{1}))){\rm \; }\mathop{ \rightarrow}\limits^{k_{1}+1 \rightarrow k_{1}}{\rm \; } ... {\rm \; }\mathop{ \rightarrow}\limits^{k_{1}+q_{m+1}-1 \rightarrow k_{1}}\]
\[\mathop{ \rightarrow}\limits^{k_{1}+q_{m+1}-1 \rightarrow k_{1}}{\rm \; } a(k_{1} ;n_{q_{m+1}}){\rm \; }\mathop{ \rightarrow}\limits^{ k_{1} \rightarrow k_{1}+q_{m+1}} {\rm \; }a(k_{1}+q_{m+1};\# _{k_{1}} (p^{2} (k_{1}))){\rm \; }\mathop{ \rightarrow}\limits^{ k_{1}+q_{m+1} \rightarrow k_{1} }\]
\begin{spacing}{-0.3}\end{spacing}
\[\mathop{ \rightarrow}\limits^{ k_{1}+q_{m+1} \rightarrow k_{1} }{\rm \; }a(k_{1} ;n_{q_{m+1}+1})\mathop{ \rightarrow}\limits^{ k_{1} \rightarrow k_{1}+q_{m+1}+1}{\rm \; }a(k_{1}+q_{m+1}+1;\# _{k_{1}} (p^{2} (k_{1}))){\rm \; }\mathop{ \rightarrow}\limits^{ k_{1}+{q_{m+1}+1} \rightarrow k_{1}+q_{m+1}}\]
\[\mathop{ \rightarrow}\limits^{ k_{1}+{q_{m+1}+1} \rightarrow k_{1}+q_{m+1}}{\rm \; }p^{2}(k_{1}(m+2)).\]

From the above scheme, we get:

1) ${\rm H} _{(m+1)^{4}}={\rm \; \{ }a(k_{1} ;n_{1} ),{\rm .}..,{\rm \; }a(k_{1} ;n_{q_{m+1}} )\}$ , where $p^{2}(k_{1})<a(k_{1} ;n_{1} )<...<a(k_{1} ;n_{q_{m+1}} )$.

2) $C_{(m+1)^{4}}=a(k_{1} ;n_{q_{m+1}+1})$, where $C_{(m+1)^{4}}>a(k_{1} ;n_{q_{m+1}})$.
\begin{spacing}{1.5}\end{spacing}
3) $p^{2}(k_{1}(m+2)) = a(k_{1}+q_{m+1};n_{q_{m+1}})$, where $k_{1}(m+2)\mathop{=}\limits^{\eqref{(6)}}k_{1}(m+1)+q_{m+1}.$

Further, 
\[\overline{{\rm H}} _{(m+1)^{4}}=\left\{ C_{(m+1)^{4}}\right\}.\]
\[ a(k_{1}+1;\# _{k_{1}} (p^{2} (k_{1})))=a(k_{1} ;n_{1})=D(p^{2}(k_{1}(m+1)))\mathop{=}\limits^{\eqref{(9)}}D(p^{2}(k_{1}(m)))=C_{m^{4}} {\rm \;}\Rightarrow\]
\begin{equation} \label{(11)} 
\Rightarrow {\rm \;}C_{m^{4}} \in {\rm H} _{(m+1)^{4}}.
\end{equation}

\eqref{(11)} makes it clear that $C_{m^{4}} \not =  C_{(m+1)^{4}}$.
 \[D_{T_{k_{1}(m+1)}}^{*}=\left\{ a \in D_{T_{k_{1}(m+1)}} {\rm \;} | {\rm \; \;} p^{2}(k_{1}(m+1))<a \le C_{(m+1)^{4}}\right\}=\]
 \[={\rm H} _{(m+1)^{4}} \cup \overline{{\rm H}} _{(m+1)^{4}}{\rm \;}\mathop{\supset}\limits^{\eqref{(11)}}{\rm \;}\overline{{\rm H}} _{m^{4}} {\rm \;} \Rightarrow {\rm \;} \overline{{\rm H}} _{m^{4}} \subset D_{T_{k_{1}(m+1)}}^{*} {\rm \;} \mathop{\Leftrightarrow }\limits^{\eqref{(10)},{\rm \;}\eqref{(9)}} {\rm \;} D_{T_{k_{1}(m)}}^{*} \subset D_{T_{k_{1}(m)}}^{*} {\rm \;} \Rightarrow\]
 \[ \Rightarrow {\rm \;} D_{T_{k_{1}(m)}}^{*} \not =  D_{T_{k_{1}(m)}}^{*}. \]
 
 As a result, a contradiction. Theorem 3.2 is proved.

\begin{center}
\textbf{\large 4. Major findings}
\end{center}
\begin{spacing}{0.7}\end{spacing}

CONCLUSION 4.1. $(\forall m \ge 2)(\exists q \in \mathbb{M}_{T})(q \in \left(m^{4}; (m+1)^{4})\right).$

CONCLUSION 4.2.  $(\forall m \ge 3) (D(p^{2} (k_{1}(m)))\in {\rm H} _{m^{4}})$.

COROLLARY 4.3. $(\forall m \ge 3) ({\rm H} _{m^{4}} \ne \varnothing)$.

CONCLUSION 4.4.  For $m\ge 3$, there are the equalities
\[\min \limits_{\scriptsize \begin{array}{l} {{\rm \; \; \;}m^2<p<(m+1)^2}  \\  {{\rm \; \; \; \; \; \; \; \; \; \; \; \; \;}p \in \mathbb{P}} \end{array}} p=\frac{D(p^{2} (k_{1}(m)))}{p(k_{1}(m) )}=\frac{\min{\rm H}_{m^{4}}}{\text{GCD}({\rm H}_{m^{4}})},\]

\noindent where $\text{GCD}({\rm H}_{m^{4}})$ -- greatest common divisor of elements of «active» set ${\rm H}_{m^{4}}$ for number $m^{4}$.

\begin{center}
\textbf{\large 5. Conclusion}
\end{center}
\begin{spacing}{0.7}\end{spacing}

The scheme for finding elements of «active» set ${\rm H} _{m^{4} }$ and «critical» element $C_{m^{4}}$ for number $m^{4}$ at each $m \ge 3$ was reviewed. Theorem on existence of the prime number between $m^{2}$ and $(m+1)^{2}$ is proved. Thus, the truth of Legendre's conjecture is established.

\begin{center}
\textbf{\large References}
\end{center}
\begin{enumerate}
\item Garipov I. About one matrix of composite numbers and her applications //arXiv preprint arXiv:2012.15745. – 2020.

\item Garipov I. About connection of one matrix of composite numbers with Legendre's conjecture //arXiv preprint arXiv:2104.06261. – 2021.
\end{enumerate}
\end{document}